\nonstopmode
\documentclass[a4paper, 10pt]{amsart}
\usepackage{latexsym}
\usepackage{fancyhdr}
\usepackage{longtable}
\usepackage{amsmath, amssymb}
\usepackage[ansinew]{inputenc}
\usepackage[all]{xy}
\usepackage{verbatim}
\usepackage{psfrag}
\usepackage{epsfig}
\usepackage{hyperref}
\usepackage{pdfsync}
\swapnumbers
\theoremstyle{plain}
\newtheorem*{lemma*}{Lemma}

\newtheorem*{theorem*}{Theorem}
\newtheorem{theorem}[section]{Theorem}
\newtheorem*{proposition*}{Proposition}

\newtheorem*{corollary*}{Corollary}
\newtheorem{corollary}[section]{Corollary}
\theoremstyle{definition}
\newtheorem*{definition*}{Definition}

\newtheorem*{example*}{Example}

\newtheorem*{remark*}{Remark}
\newtheorem*{remarks*}{Remarks}

\newenvironment{demo}[1]{\par\smallskip\noindent{\bf #1.}}{\par\smallskip}

\def\al{\alpha}

\def\la{\lambda}

\def\rh{\rho}

\def\si{\sigma}

\def\Si{\Sigma}

\def\R{\mathbb{R}}

\def\<{\langle}
\def\>{\rangle}
\def\-{\backslash}
\renewcommand{\o}{\circ}

\let\on=\operatorname
\def\sr#1%
{\ifmmode{}^\dagger\else${}^\dagger$\fi\ifvmode
\vbox to 0pt{\vss
 \hbox to 0pt{\hskip\hsize\hskip1em
 \vbox{\hsize3cm\eightpoint\raggedright\pretolerance10000
 \noindent #1\hfill}\hss}\vss}\else
 \vadjust{\vbox to0pt{\vss%
 \hbox to 0pt{\hskip\hsize\hskip1em%
 \vbox{\hsize3cm\eightpoint\raggedright\pretolerance10000%
 \noindent #1\hfill}\hss}\vss}}\fi%
}

\def\ignore#1{}

\def\mapsfrom{\DOTSB\kern.2em%
\setbox0=\hbox{$\leftarrow$\kern-.15em\raise0.1ex\hbox{$\shortmid$}}\box0%
\kern.3em}

\newcommand{\nmb}[2]{\ifx!#1{\ref{nmb:#2}}\else{\label{nmb:#2}}\fi}
\allowdisplaybreaks
\date{April 29, 2010}

\title[Addendum: ``Lifting smooth curves over invariants, III'']
{Addendum to: ``Lifting smooth curves over invariants for representations of compact Lie groups, III'' [J. Lie Theory \textbf{16} (2006), No. 3, 579--600.]}

\author
[A.~Kriegl, M.~Losik, P.W.~Michor, A.~Rainer]
{Andreas Kriegl, Mark Losik, Peter W. Michor, and Armin Rainer}

\address{Andreas Kriegl: Fakult\"at f\"ur Mathematik, Universit\"at Wien, 
Nordbergstrasse~15, A-1090 Wien, Austria}
\email{andreas.kriegl@univie.ac.at}

\address{Mark Losik: Saratov State University, 
ul. Astrakhanskaya, 83, 410026 Saratov, Russia}
\email{losikMV@info.sgu.ru}

\address{Peter W.\ Michor: Fakult\"at f\"ur Mathematik, Universit\"at Wien, 
Nordbergstrasse~15, A-1090 Wien, Austria}
\email{peter.michor@univie.ac.at}

\address{Armin Rainer: Fakult\"at f\"ur Mathematik, Universit\"at Wien,
Nordbergstrasse 15, A-1090 Wien, Austria}
\email{armin.rainer@univie.ac.at}

\thanks{AK was supported by FWF-Project P~23082-N13.
PM was supported by FWF-Project P~21030-N13.
AR was supported by FWF-Projects J~2771-N13 and P~22218-N13.}

\begin{document}

\begin{abstract}
We improve the main results in the paper from the title using a recent refinement of
Bronshtein's theorem 
due to Colombini, Orr\'u, and Pernazza.
They are then in general best possible both in the hypothesis and in the
outcome. As a consequence we obtain a result on lifting smooth mappings 
in several variables.
\end{abstract}

\maketitle

A recent refinement of Bronshtein's theorem \cite{Bronshtein79} and of some of its consequences due to
Colombini, Orr\'u, and Pernazza \cite{ColombiniOrruPernazza08}
(namely theorem \ref{roots}(i) below) allows to essentially improve our
main results in \cite{KLMR06}; see theorem \ref{liftsf}
and
corollary \ref{liftsp} below. The improvement consists in weakening the
hypothesis considerably: In \cite{KLMR06} we needed a curve $c$ to be of
class
\begin{enumerate}
 \item[(i)] $C^k$ in order to admit a differentiable lift with locally
 bounded derivative,
 \item[(ii)] $C^{k+d}$ in order to admit a $C^1$-lift, and
 \item[(iii)] $C^{k+2d}$ in order to admit a twice differentiable lift.
\end{enumerate}
It turns out that theorem \ref{liftsf} and corollary \ref{liftsp} are in
general best possible both in the hypothesis and in the outcome.
In theorem \ref{thm:4} and corollary \ref{cor:5} we deduce some results on
lifting smooth mappings in several variables.

\subsection*{Refinement of Bronshtein's theorem}
Bronshtein's theorem \cite{Bronshtein79} (see also Wakabayashi's version
\cite{Wakabayashi86}) states that,
for a curve of monic hyperbolic polynomials
\begin{equation} \label{P}
P(t)(x) = x^n + \sum_{j=1}^n (-1)^j a_j(t) x^{n-j}.
\end{equation}
with coefficients $a_j \in C^n(\R)$ ($1 \le j \le n$), there exist
differentiable functions $\la_j$ ($1 \le j \le n$) with locally bounded
derivatives
which parameterize the roots of $P$. A polynomial is called hyperbolic if
all its roots are real.

The following theorem refines Bronshtein's theorem \cite{Bronshtein79} and also a result of
Mandai \cite{Mandai85} and a result of Kriegl, Losik, and Michor
\cite{KLM04}.
In \cite{Mandai85} the coefficients are required to be of class $C^{2n}$ for
$C^1$-roots, and in \cite{KLM04} they are assumed to be $C^{3n}$ for twice
differentiable roots.

\begin{theorem}[{\cite[2.1]{ColombiniOrruPernazza08}}] \label{roots}
Consider a curve $P$ of monic hyperbolic polynomials \eqref{P}. Then:
\begin{enumerate}
 \item[(i)] If $a_j \in C^n(\R)$ ($1 \le j \le n$), then there exist
 functions $\la_j \in C^1(\R)$ ($1 \le j \le n$) which parameterize the
 roots of $P$.
 \item[(ii)] If $a_j \in C^{2n}(\R)$ ($1 \le j \le n$), then the roots of
 $P$ may be chosen twice differentiable.
\end{enumerate}
\end{theorem}

Counterexamples (e.g. in \cite[section 4]{ColombiniOrruPernazza08}) show that in
this result the assumptions on $P$ cannot be weakened.

\subsection*{Improvement of the results in \cite{KLMR06}}
Let $\rh : G \to \on{O}(V)$ be an orthogonal representation of a compact
Lie group $G$ in a real finite dimensional Euclidean vector space $V$.
Choose a minimal system of
homogeneous generators $\si_1,\ldots,\si_n$ of the algebra
$\mathbb{R}[V]^G$ of $G$-invariant polynomials on V.
Define
\[
d=d(\rh) := \max \{\deg \si_i : 1 \le i \le n \},
\]
which is independent of the choice of the $\si_i$ (see \cite[2.4]{KLMR06}).

If $G$ is a finite group, we
write $V = V_1 \oplus \cdots \oplus V_l$ as orthogonal direct sum of
irreducible subspaces $V_i$.
We choose $v_i \in V_i \backslash \{0\}$ such that the cardinality of the
corresponding isotropy group $G_{v_i}$ is maximal,
and put
\[
k = k(\rh) := \max\{d(\rho), |G|/|G_{v_i}| : 1 \le i \le l\}.
\]

The mapping $\si = (\si_1,\ldots,\si_n) : V \to \R^n$ induces a
homeomorphism between the orbit space $V/G$ and the image $\si(V)$.
Let $c : \R \to V/G = \si(V) \subseteq \R^n$ be a smooth curve in the orbit
space (smooth as a curve in $\R^n$).
A curve $\bar c : \R \to V$ is called lift of $c$ if $\si \o \bar c = c$.
The problem of lifting curves smoothly over invariants is independent of the
choice of the $\si_i$ (see \cite[2.2]{KLMR06}).

\begin{theorem} \label{liftsf}
Let $\rh : G \to \on{O}(V)$ be a representation of a finite group $G$. Let
$d=d(\rh)$ and $k=k(\rh)$.
Consider a curve $c : \R \to V/G = \si(V) \subseteq \R^n$ in the orbit
space of $\rh$.
Then:
\begin{enumerate}
 \item[(i)] If $c$ is of class $C^k$, then any differentiable lift $\bar c
 : \R \to V$ of $c$ (which always exists) is actually $C^1$.
 \item[(ii)] If $c$ is of class $C^{k+d}$, then there exists a global twice
 differentiable lift $\bar c : \R \to V$ of $c$.
\end{enumerate}
\end{theorem}

\begin{demo}{Proof}
  (i) Let $\bar c$ be any differentiable lift of $c$.
  Note that the existence of $\bar c$ is guaranteed for any $C^d$-curve $c$,
  by \cite{KLMR05}.
  In the proof of \cite[8.1]{KLMR06} we construct curves of monic hyperbolic
  polynomials $t \mapsto P_i(t)$ which have the regularity of $c$ and whose
  roots
  are parameterized by $t \mapsto \< v_i \mid g.\bar c(t) \>$ ($g \in
  G_{v_i}\backslash G$).

  If $c$ is of class $C^k$, then theorem \ref{roots}(i) provides $C^1$-roots
  of $t \mapsto P_i(t)$. By the proof of \cite[4.2]{KLMR06} we obtain that
  the parameterization $t \mapsto \< v_i \mid g.\bar c(t) \>$ is $C^1$ as
  well. Hence $\bar c$ is a $C^1$-lift of $c$.
  Alternatively, the proof of \ref{roots}(i) in
  \cite{ColombiniOrruPernazza08} actually shows that any differentiable
  choice of roots is $C^1$.

  (ii) Let $c$ be of class $C^{k+d}$. The existence of a global twice
  differentiable lift $\bar c$ of $c$ follows from the proof of \cite[5.1 and
  5.2]{KLMR06},
  where we use (i) instead of \cite[4.2]{KLMR06}.
\qed\end{demo}

\begin{corollary} \label{liftsp}
Let $\rh : G \to \on{O}(V)$ be a polar representation of a compact Lie
group $G$.
Let $\Si \subseteq V$ be a section, $W(\Si)=N_G(\Si)/Z_G(\Si)$ its
generalized Weyl group, and $\rh_\Si : W(\Si) \to \on{O}(\Si)$ the induced
representation.
Let $d=d(\rh_\Si)$ and $k=k(\rh_\Si)$.
Consider a curve $c : \R \to V/G = \si(V) \subseteq \R^n$ in the orbit
space of $\rh$.
Then:
\begin{enumerate}
 \item[(i)] If $c$ is of class $C^k$, then there exists a global orthogonal
 $C^1$-lift $\bar c : \R \to V$ of $c$.
 \item[(ii)] If $c$ is of class $C^{k+d}$, then there exists a global
 orthogonal twice differentiable lift $\bar c : \R \to V$ of $c$. \qed
\end{enumerate}
\end{corollary}

The examples which show that the hypothesis in \ref{roots} are best
possible also imply that in general the hypothesis in \ref{liftsf} and
\ref{liftsp} cannot be improved.

On the other hand the outcome of \ref{liftsf} and \ref{liftsp} cannot be
refined either: A $C^\infty$- curve $c$ does in general not allow a
$C^{1,\al}$-lift for any $\al>0$. See \cite{Glaeser63R}, \cite{AKLM98},
\cite{BBCP06}.
But see  also \cite{Bony05} and \cite[remark 4.2]{KLMR06}.

Note that the improvement affects also \cite[part 6]{LR07}.

\subsection*{Lifting smooth mappings in several variables}
From theorem \ref{liftsf} we can deduce a lifting result for mappings in several variables.

\begin{theorem}\label{thm:4}
  Let $\rh : G \to \on{O}(V)$ be a representation of a finite group $G$, 
  $d=d(\rh)$, and $k=k(\rh)$. Let $U \subseteq \R^q$ be open.
  Consider a mapping $f : U \to V/G = \si(V) \subseteq \R^n$ of class $C^k$.
  Then any continuous lift $\bar f : U \to V$ of $f$ is actually locally Lipschitz.
\end{theorem}

\begin{demo}{Proof}
  Let $c : \R \to U$ be a $C^\infty$-curve. By theorem \ref{liftsf}(i) 
  the curve $f \o c$ admits a $C^1$-lift $\overline{f \o c}$. A further continuous lift of 
  $f \o c$ is formed by $\bar f \o c$. 
  By \cite[5.3]{LMRac} we can conclude that $\bar f \o c$ is locally Lipschitz.
  So we have shown that $\bar f$ is locally Lipschitz along $C^\infty$-curves. 
  By Boman \cite{Boman67} (see also \cite[12.7]{KM97}) that implies that $\bar f$ is locally Lipschitz. 
\qed\end{demo}

In general there will not always exist a continuous lift of $f$ (for instance, if $G$ is a finite rotation group and $f$ is defined near $0$).
However, if $G$ is a finite reflection group, then any continuous $f$ allows a continuous lift 
(since the orbit space can be embedded homeomorphically in $V$).

\begin{corollary}\label{cor:5}
  Let $\rh : G \to \on{O}(V)$ be a polar representation of a compact connected Lie
  group $G$.
  Let $\Si \subseteq V$ be a section, $W(\Si)=N_G(\Si)/Z_G(\Si)$ its
  generalized Weyl group, $\rh_\Si : W(\Si) \to \on{O}(\Si)$ the induced
  representation, $d=d(\rh_\Si)$, and $k=k(\rh_\Si)$. Let $U \subseteq \R^q$ be open.
  Consider a mapping $f : U \to V/G = \si(V) \subseteq \R^n$ of class $C^k$.
  Then there exists an orthogonal lift $\bar f : U \to V$ of $f$ which is locally Lipschitz.
\end{corollary}

\begin{demo}{Proof}
  The Weyl group $W(\Si)$ is a finite reflection group, since $G$ is connected. 
\qed\end{demo}

\def\cprime{$'$}
\providecommand{\bysame}{\leavevmode\hbox to3em{\hrulefill}\thinspace}
\providecommand{\MR}{\relax\ifhmode\unskip\space\fi MR }
\providecommand{\MRhref}[2]{%
  \href{http://www.ams.org/mathscinet-getitem?mr=#1}{#2}
}
\providecommand{\href}[2]{#2}


\end{document}